\documentclass[a4paper,reqno, 12pt,titlepage]{amsart}
\usepackage{amssymb, graphicx, enumerate, enumitem, color}
\usepackage[foot]{amsaddr}

\usepackage{tikz}
\usetikzlibrary{calc}
\usepackage[pagewise]{lineno}

\usepackage[unicode=true]{hyperref}
\hypersetup{
	pdftitle={Separating Tree-chromatic number from 
		Path-chromatic Number},
	pdfauthor={Fidel Barrera-Cruz, Stefan Felsner, Tam\'{a}s M\'{e}sz\'{a}ros, Piotr Micek, Heather Smith, Libby Taylor, William T. Trotter},
	colorlinks
}

\newtheorem{theorem}{Theorem}[section]

\theoremstyle{remark}

\newcommand{\cgB}{\mathcal{B}}

\newcommand{\pt}{\operatorname{pt}}

\newcommand{\Ram}{\operatorname{Ram}}

\newcommand{\DMR}{\operatorname{DMR}}
\newcommand{\DML}{\operatorname{DML}}
\newcommand{\OMR}{\operatorname{OMR}}
\newcommand{\OML}{\operatorname{OML}}
\newcommand{\IFS}{\operatorname{IFS}}
\newcommand{\ISF}{\operatorname{ISF}}

\DeclareMathOperator\pchr{path-\chi}
\DeclareMathOperator\tchr{tree-\chi}

\let\le\leqslant
\let\ge\geqslant
\let\leq\leqslant
\let\geq\geqslant

\let\epsilon\varepsilon

\renewenvironment{enumerate}{\begin{enumorig}[label=\textup{(\roman*)}, noitemsep, topsep=2pt plus 2pt, labelindent=.2em, leftmargin=*, widest=iii]}{\end{enumorig}}

\begin{document}


\title[SEPARATING TREE-$\chi$ FROM PATH-$\chi$]%
{Separating Tree-chromatic number from 
Path-chromatic Number}

\author[BARRERA-CRUZ]{Fidel Barrera-Cruz}
\address[F.~Barrera-Cruz, H.~Smith, L.~Taylor, W.~T.~Trotter]{School of Mathematics\\
  Georgia Institute of Technology\\
  Atlanta, Georgia 30332}
\email{fidelbc@math.gatech.edu}

\author[FELSNER]{Stefan Felsner}
\address[S.~Felsner]{Institut f\"ur Mathematik\\
  Technische Universit\"{a}t Berlin\\
  Strasse des 17. Juni 136\\
  D-10623 Berlin\\
  Germany}
\email{felsner@math.tu-berlin.de}

\author[M\'{E}SZ\'{A}ROS]{Tam\'{a}s M\'{e}sz\'{a}ros}
\address[T.~M\'{e}sz\'{a}ros]{Fachbereich Mathematik und Informatik\\
  Kombinatorik und Graphentheorie\\
  Freie Universit\"at Berlin\\
  Arnimallee 3, 14195 Berlin\\
  Germany}
\email{meszaros.tamas@fu-berlin.de}

\author[MICEK]{Piotr Micek}
\address[P.~Micek]{Theoretical Computer Science Department\\
  Faculty of Mathematics and Computer Science\\
  Jagiellonian University\\
  Krak\'ow, Poland}
\email{piotr.micek@tcs.uj.edu.pl}

\author[SMITH]{Heather Smith}
\email{heather.smith@math.gatech.edu}

\author[TAYLOR]{Libby Taylor}
\email{libbytaylor@gatech.edu}

\author[TROTTER]{William T. Trotter}
\email{trotter@math.gatech.edu}

\date{January 30, 2019}

\keywords{Binary tree, Ramsey theory, tree-decomposition}

\thanks{Tam\'{a}s M\'{e}sz\'{a}ros is supported by the DRS 
Fellowship Program at Freie Universit\"{a}t Berlin.}
\thanks{Piotr Micek is partially supported by the National Science Center of Poland under grant no.\ 2015/18/E/ST6/00299.}
\thanks{Stefan Felsner is partially supported by DFG Grant Fe 340/11-1.}
 
\begin{abstract}
We apply Ramsey theoretic tools to 
show that there is a family of graphs which have tree-chromatic 
number at most~$2$ while the path-chromatic number is unbounded.
This resolves a problem posed by Seymour.
\end{abstract}

\maketitle

\section{Introduction}

Let $G$ be a graph.
A \textit{tree-decomposition} of $G$ is a pair $(T,\cgB)$ where $T$ 
is a tree and $\cgB=(B_t\mid t\in V(T))$ is a family of subsets of 
$V(G)$, satisfying:
\begin{enumerate}
\item[(T1)] for each $v\in V(G)$ there exists $t\in V(T)$ with $v\in B_t$; 
and for every edge $uv \in E(G)$ there exists $t\in V(T)$ with $u,v\in B_t$;
\item[(T2)] for each $v\in V(G)$, if $v\in B_t \cap B_{t''}$ for some 
$t,t''\in V(T)$, and $t'$ lies on the path in $T$ between $t$ and $t''$, 
then $v\in B_{t'}$.
\end{enumerate}
Many researchers refer to the subset $B_t$ as a \emph{bag} and they
consider $B_t$ as an induced subgraph of $G$.  With this convention,
$|B_t|$ is just the number of vertices of $G$ in the bag $B_t$,
while $\chi(B_t)$ is the chromatic number of the induced subgraph
of $G$ determined by the vertices in $B_t$.

The quality of a tree-decomposition $(T,(B_t\mid t\in V(T)))$ is usually 
measured by its \emph{width}, i.e.\ the maximum of $|B_t|-1$ over all $t\in V(T)$.
Then the \emph{tree-width} of $G$ is the minimum width of a tree-decomposition of $G$.
In this paper we study the tree-chromatic number of a graph,
a concept introduced by Seymour in~\cite{bib:Seym}. The \emph{chromatic number} of 
a tree-decomposition $(T,(B_t\mid t\in V(T)))$ is the maximum of 
$\chi(B_t)$ over all $t\in V(T)$.  The \emph{tree-chromatic number} of 
$G$, denoted by $\tchr(G)$, is the minimum chromatic number of a 
tree-decomposition of $G$. A tree-decomposition $(T,(B_t\mid t\in V(T)))$ 
is a \emph{path-decomposition} when $T$ is a path.
The \emph{path-chromatic number} of $G$, denoted by $\pchr(G)$, is 
the minimum chromatic number of a path-decomposition of $G$.
Clearly, for every graph $G$ we have
\[
\omega(G)\le\tchr(G)\le\pchr(G)\le\chi(G).
\]

Furthermore, if $G=K_n$ is the complete graph on $n$ vertices,
then $\omega(G)=\chi(G)=n$, so all these inequalities can be
tight.  Accordingly, it is of interest to ask whether for consecutive 
parameters in this series of inequalities, there is a sequence of graphs 
for which one parameter is bounded while the next parameter is
unbounded.

In~\cite{bib:Seym}, Seymour proved that the classic Erd\H{o}s 
construction~\cite{bib:Erdo} for graphs with large girth and large 
chromatic number yields a sequence $\{G_n:n\ge1\}$ with $\omega(G_n)=2$ 
and $\tchr(G_n)$ unbounded.

For an integer $n\ge2$, the \emph{shift graph} $S_n$ is a graph whose
vertex set consists of all closed intervals of the form
$[a,b]$ where $a$, $b$ are integers with $1\le a<b\le n$.
Vertices $[a,b]$, $[c,d]$ are adjacent in $S_n$ when $b=c$ or $d=a$.
As is well known (and first shown in~\cite{bib:ErdHaj}), 
$\chi(S_n)=\lceil\lg n\rceil$, for every $n\ge2$.  
On the other hand, $S_n$ has a natural simple path decomposition. $T$ is simply the path on the vertices $t_1,t_2\dots, t_n$ and for $1\leq i\leq n$ we have $B_{t_i}=\{[a,b]\in V(S_n)\ :\ a\leq i\leq b\}$. Then for every $1\leq i\leq n$ the bag $B_{t_i}$ is the union of two independent sets, namely $\{[a,b]\in V(S_n)\ :\ a<i\leq b\}$ and $\{[a,b]\in V(S_n)\ :\ a\leq i< b\}$, and hence the chromatic number of the corresponding induced subgraph is at most $2$. This shows that $\pchr(S_n)\leq2$ for every $n\geq2$, so as noted in~\cite{bib:Seym},
the family of shift graphs has bounded path-chromatic number and 
unbounded chromatic number.

Accordingly, it remains only to determine whether there is
an infinite sequence of graphs with bounded tree-chromatic number
and unbounded path-chromatic number.  However, these two parameters
appear to be more subtle in nature.  As a first step,
Huynh and Kim~\cite{bib:HuyKim} showed that there is
an infinite sequence $\{G_n:n\ge1\}$ of graphs with
$\tchr(G_n)\rightarrow\infty$ and $\tchr(G_n)<\pchr(G_n)$ for all 
$n\ge1$.

In~\cite{bib:Seym}, Seymour proposed the following
construction.   Let $T_n$ be the complete (rooted) binary tree with
$2^n$ leaves.  When $y$ and $z$ are distinct vertices in
$T_n$, the path from $y$ to $z$ is called
a ``$V$'' when the unique point
on the path which is closest to the root of $T_n$ is an intermediate
point $x$ on the path which is \emph{strictly} between $y$ and $z$.
We refer to $x$ as the \emph{low point} of the $V$ formed by 
$y$ and $z$.

For a fixed value of $n$, we can then form a graph $G_n$ whose
vertices are the $V$'s in $T_n$.  We take $V$ adjacent to
$V'$ in $G_n$ when an end point of one of the two paths is the low point
of the other.  Clearly, $\omega(G_n)\le 2$. Furthermore, it is easy to see that 
$\chi(G_n)\rightarrow \infty$ with $n$ (we will say more about 
this observation later in the paper), and Seymour~\cite{bib:Seym} suggested 
that the family $\{G_n:n\ge1\}$ has unbounded path-chromatic number. 

However, we will show that graphs in the family $\{G_n:n\ge1\}$ have
bounded path-chromatic number.  In fact, we will use Ramsey theoretic
tools developed by Milliken~\cite{bib:Mill} to show that if we fix
$r\ge2$, and assume we have a path-decomposition of $G_n$ witnessing
that $\pchr(G_n)\le r$, then this decomposition is (essentially) uniquely 
determined.  Furthermore, this decomposition actually witnesses that 
$\pchr(G_n)\le2$.

Moreover, in analyzing this decomposition, we discovered the following
minor modification.  In the binary tree $T_n$,
a subtree is called a ``$Y$'' when it has~$3$ leaves and the closest vertex in
the subtree to the root of $T_n$ is one of the three leaves.  
We then let $H_n$ be the graph whose vertex set consists of the
$V$'s and $Y$'s in $T_n$.  Furthermore,  $Y$ is adjacent to $Y'$ in
$H_n$ if and only if one of the two upper leaves of one of them is the lowest
leaf in the other.  Also, a $Y$ is adjacent to a $V$ if and only if one of
the upper leaves in the $Y$ is the low point of the $V$.

It is clear from the natural tree-decomposition of $H_n$
that $\tchr(H_n)\le 2$.  Using Ramsey theoretic tools, we will
then be able to show that $\pchr(H_n)\rightarrow\infty$ with $n$, so
that Seymour's question has been successfully resolved.

\section{Ramsey Theory on Binary Trees}\label{sec:ramsey}
 
The Ramsey theoretic concepts discussed here are treated in a more
comprehensive manner by Milliken~\cite{bib:Mill}\footnote{The particular
result we need is Theorem~2.1 on page 220.  Note that Milliken
credits the result to Halpern, L\"{a}uchi, Laver and Pincus and
comments on the history of the result.}, but we will find
it convenient to use somewhat different notation and terminology.

For a positive integer $n$, we view the \textit{complete binary tree 
$T_n$} as the poset\footnote{The complete (rooted) binary tree we 
discussed in an informal manner in the opening section of this paper 
is just the cover graph of the poset $T_n$ defined here.}
consisting of all binary strings of length at most~$n$, 
with $x\le y$ in $T_n$ when $x$ is an initial segment in $y$.  
The empty string, denoted $\emptyset$, 
is then the zero (least element) of $T_n$.  For all $n\ge1$, $T_n$ has
$2^{n+1}-1$ elements and height~$n+1$.
In particular, $T_0$ is the one-point
poset whose only element is the empty string.

When $n\ge1$ and $x$ is a binary string of length~$n$,
we will denote coordinate $i$ of $x$ as $x(i)$ and
when a string is of modest length, we may write it explictly, e.g.,
$x=01001101$.  When $n\ge m>p\ge0$, $x$ is a string of length $p$,
$y$ is a string of length~$m$ and $x<y$ in $T_n$, we say $y$ is
in the \textit{left tree above $x$} when $y(p+1)=0$ and we say
$y$ is in the \textit{right tree above $x$} when $y(p+1)=1$.

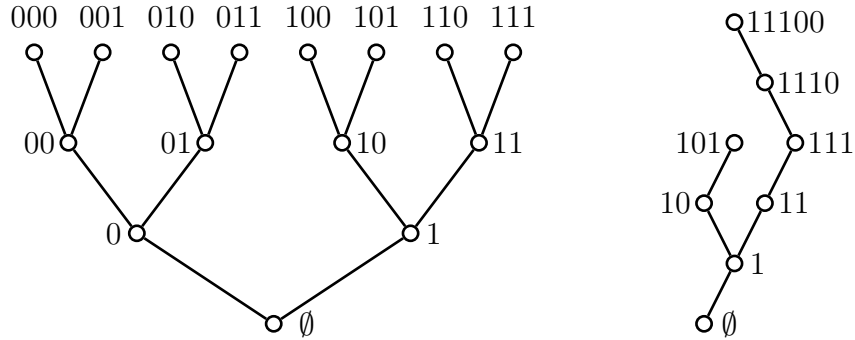
\begin{figure}
\begin{center}
\begin{center}
\begin{tikzpicture}
[line width=1pt,level distance=12mm,
every node/.style={draw, color=black, line width=1pt, fill=white, circle, minimum size=2mm, inner sep=0cm, label distance=0mm},
level 1/.style={sibling distance=36mm},
level 2/.style={sibling distance=18mm},
level 3/.style={sibling distance=9mm}]
\node [label={[label distance=1mm]0:$\emptyset$}] {} [grow=up]
child {node [label={0:$1$}] {}
	child {node [label={0:$11$}] {}
		child {node [label={90:$111$}] {}
		}
		child {node [label={90:$110$}] {}
		}
	}
	child {node [label={0:$10$}] {}
		child {node [label={90:$101$}] {}
		}
		child {node [label={90:$100$}] {}
		}
	}
}
child {node [label={180:$0$}] {}
	child {node [label={180:$01$}] {}
		child {node [label={90:$011$}] {}
		}
		child {node [label={90:$010$}] {}
		}
	}
	child {node [label={180:$00$}] {}
		child {node [label={90:$001$}] {}
		}
		child {node [label={90:$000$}] {}
		}
	}
}
;
\end{tikzpicture}
\hspace{.5in}
\begin{tikzpicture}
[line width=1pt,level distance=8mm,
every node/.style={draw, color=black, line width=1pt, fill=white, circle, minimum size=2mm, inner sep=0cm, label distance=0mm},
level 1/.style={sibling distance=8mm},
level 2/.style={sibling distance=8mm},
level 3/.style={sibling distance=8mm},
level 4/.style={sibling distance=8mm}]
\node [label={0:$\emptyset$}] {} [grow=up]
child {node [label={0:$1$}] {}
	child {node [label={0:$11$}] {}
		child {node [label={0:$111$}] {}
			child [missing]
			child {node [label={0:$1110$}] {}
				child [missing]
				child {node [label={0:$11100$}] {}
				}
			}
		}
		child [missing]
	}
	child {node [label={180:$10$}] {}
		child {node [label={180:$101$}] {}
		}
		child [missing]
	}
}
child [missing]
;
\end{tikzpicture}
\end{center}
\end{center}
\caption{Binary Trees: Down sets in $T_n$\label{fig:binary-trees}}
\end{figure}

Recall that in a poset $P$, a subposet $Q$ of $P$ is called a \textit{down
set} if $x\in Q$ whenever $y\in Q$ and $x\le y$ in $P$. We will refer
to down sets of the complete binary tree $T_n$ as \textit{binary trees}.
In Figure~\ref{fig:binary-trees}, we show on the left the 
complete binary tree $T_3$.  On the right, we show a binary tree $Q$ which
will be a down set in any complete binary tree $T_n$ with $n\ge5$.
 
Let $n\ge0$, let $Q$ be a binary tree in $T_n$, and let
$R$ be a subposet of $T_n$. Following Milliken~\cite{bib:Mill}, we will say
$R$ is a \textit{strong copy of $Q$} when
there is a function $f:Q\rightarrow R$ satisfying the following
two requirements:

\begin{enumerate}
\item $f$ is a poset isomorphism, i.e., $f$ is a bijection
and for all $x,y\in Q$, $x\le y$ in $Q$ if and only if
$f(x)\le f(y)$ in $R$.
\item For all $x,y\in Q$ with $x<y$ in $Q$, $y$ is in the
left tree above $x$ in $Q$ if and only if $f(y)$ is in the
left tree above $f(x)$ in $T_n$.
\end{enumerate}

Since we are concerned with binary trees, we note that when $f$ satisfies
the preceding two conditions, then it automatically implies that $y$ is
in the right tree above $x$ if and only if $f(y)$ is in the right tree
above $f(x)$.  

For the remainder of this paper, when $r\ge1$, we let
$[r]$ denote the set $\{1,2,\dots,r\}$.  Also, an $r$-coloring
of a set $X$ is just a map $\Phi:X\rightarrow[r]$.  In some situations,
we will consider a coloring $\Phi$ using a set of $r$ colors, but
the set will not simply be the set $[r]$.

The following result is a straightforward extension of
the special case of Theorem~2.1 from~\cite{bib:Mill} for binary
trees.

\begin{theorem}\label{thm:ramsey-bt}
For every triple $(Q,p,r)$, where $Q$ is a binary tree, and $p$ and $r$ are positive integers
with $p$ at least as large as the height of $Q$,
there is a least positive integer $n_0=\Ram(Q,p,r)$ 
so that if $n\ge n_0$ and $\Phi$ is an $r$-coloring of the 
strong copies of $Q$ in $T_n$, then there is a color $\alpha\in[r]$
and a subposet $R$ of $T_n$ such that $R$ is a strong copy of
$T_p$ and $\Phi$ assigns color $\alpha$ to every strong copy of
$Q$ contained in $R$.
\end{theorem}

\section{Separating Tree-chromatic Number and Path-chromatic Number}

For the remainder of the paper, for a positive integer $n$,
we let $G_n$ be the graph of the $V$'s in the complete binary
tree $T_n$. Strictly speaking, a 
vertex $V$ in $G_n$ is a path which is determined
by its two endpoints, but we find it convenient to
specify $V$ as a triple $(x,y,z)$, where $y$ and
$z$ are the endpoints of the path and $x$ is the low point
on the path.  We view $V$ as a triple and not a $3$-element set so
we can follow the convention that $y$ is in
the left tree above $x$ and $z$ is in the right tree above $x$.
When $V_1=(x_1,y_1,z_1)$ and $V_2=(x_2,y_2,z_2)$ are vertices
in $G_n$, we note that $V_1$ and $V_2$ are adjacent if and only if one of the
following four statements holds: $z_1=x_2$, $y_1=x_2$, $y_2=x_1$ or 
$z_2=x_1$.

Also, for each $n\ge1$, we let $H_n$ be the graph of $V$'s and $Y$'s in $T_n$.
Of course, $G_n$ is an induced subgraph of $H_n$.  Furthermore, the 
natural tree-decomposition of $H_n$ shows that $\tchr(H_n)\le 2$ for all
$n\ge1$.  

Our goals for this section are to prove the following two theorems.

\begin{theorem}\label{thm:pchrGn}
For all $n\ge1$, the path-chromatic number of the graph $G_n$ of
$V$'s in the complete binary tree $T_n$ is at most~$2$.
\end{theorem}

\begin{theorem}\label{thm:pchrHn}
For every positive integer $r$, there is a least positive integer
$n_0$ so that if $n\ge n_0$,  then the path-chromatic number of
the graph $H_n$ of $V$'s and $Y$'s in the complete binary tree $T_n$
has chromatic number larger than $r$.
\end{theorem}

We elect to follow the line of our research and prove the second of these
two theorems first.  In accomplishing
this goal, we will discover a path-decomposition of $G_n$ witnessing
that $\pchr(G_n)\le 2$ for all $n\ge1$.  

Our argument for Theorem~\ref{thm:pchrHn} will proceed by contradiction, i.e.
we will assume that there is some positive integer $r$ such
that $\pchr(H_n)\le r$ for all $n\ge1$.  The contradiction will
come when $n$ is sufficiently large in comparison to $r$.  

For the moment, we take $n$ as a large but unspecified integer.
Later, it will be clear how large $n$ needs to be.
We then take a path-decomposition of $H_n$ witnessing that $\pchr(H_n)\le r$.
We may assume that the host path in this decomposition is the set $\mathbb{N}$ 
of positive integers with $i$ adjacent to $i+1$ in $\mathbb{N}$
for all $i\ge1$.  For each vertex $v$ in $H_n$, the
set of all integers $i$ for which $v\in B_i$ is a set
of consecutive integers, and we denote the least integer in this
set as $a_v$ and the greatest integer as $b_v$.  Abusing notation
slightly, we will denote this set as $[a_v,b_v]$, i.e., this
interval notation identifies the integers $i\in \mathbb{N}$ with $a_v\le i\le b_v$.
Alternatively, $[a_v,b_v]$ is just the set of integers $i$ for which
$v$ is in the bag $B_i$. We point out the requirement that 
$[a_v,b_v]\cap[a_u,b_u] \neq\emptyset$ when $v$ and $u$ are adjacent 
vertices in $H_n$.

After duplicating vertices in the path decomposition if necessary, we may
assume that $a_v < b_v$ for every vertex
$v\in V(H_n)$.  Similarly, after adding extra vertices to the path decomposition 
if necessary, we may assume that for each integer $i$, there
is at most one vertex $v\in V(H_n)$ with $i\in\{a_v,b_v\}$. 

For each $i\in \mathbb{N}$, we let $G_n(i)$ denote the
induced subgraph of $G_n$ determined by those vertices $v\in G_n$ with 
$i\in[a_v,b_v]$.  Alternatively, $G_n(i)$ is the subgraph of $G_n$ induced
by the vertices in bag $B_i$.  The graph $H_n(i)$ is defined analogously.  

We pause here to point out an essential detail for the remainder of the proof.
Since $\chi(G_n(i))\le \chi(H_n(i))\le r$ for all integers $i$, then
for all $q>2^r$, there is no positive integer $i$ for which either
$G_n(i)$ or $H_n(i)$ contains the shift graph $S_q$ as a subgraph.  

To begin to make the connection with Ramsey theory, we observe
that there is a natural $1$--$1$ correspondence between
$V$'s in $G_n$ and strong copies of $T_1$ in $T_n$.
So in the discussion to follow, we will interchangeably view 
a vertex $V=(x,y,z)$ of $G_n$ as a path in $T_n$ and as a 
$3$-element subposet of $T_n$ forming a strong copy of $T_1$.
Of course, we are abusing notation slightly by referring to $T_n$
as a graph and as a poset, but by now the notion that as a graph,
we are referring to the cover graph of the poset should be clear.

In the discussion to follow, when we discuss a family
$\{V_j:j\in[t]\}$ of $V$'s in $G_n$, we will let
$V_j=(x_j,y_j,z_j)$, and we will let $[a_j,b_j]$ be the
interval in the path-decomposition corresponding to $V_j$, for
each $j\in[t]$.

Let $(V_1,V_2)$ be an ordered pair of vertices in $G_n$.
Referring to the binary trees in Figure~\ref{fig:Q's}, we consider the $7$ 
different ways this pair can appear in $T_n$ so that the two paths
$V_1$ and $V_2$ have at most one vertex from $T_n$ in common:

\begin{figure}
\begin{center}
\begin{center}
\begin{tikzpicture}
[line width=1pt,level distance=8mm,
every node/.style={draw, color=black, line width=1pt, fill=white, circle, minimum size=2mm, inner sep=0cm, label distance=0mm},
level 1/.style={sibling distance=8mm},
level 2/.style={sibling distance=8mm},
level 3/.style={sibling distance=8mm},
level 4/.style={sibling distance=8mm}]
\node  (r0) [label={0:$\emptyset$}] {}  [grow=up]
child {node [label={0:$1$}] {}
	child {node [label={0:$11$}] {}
	}
	child {node [label={180:$10$}] {}
	}
}
child {node [label={180:$0$}] {}
}
;
\node at (r0) [label={[label distance = 1mm]270:$Q_{1}$}] {};
\end{tikzpicture}
\hspace{0.5in}
\begin{tikzpicture}
[line width=1pt,level distance=8mm,
every node/.style={draw, color=black, line width=1pt,fill=white,circle, minimum size=2mm,  inner sep=0cm, label distance=0mm},
level 1/.style={sibling distance=8mm},
level 2/.style={sibling distance=8mm},
level 3/.style={sibling distance=8mm},
level 4/.style={sibling distance=8mm}]
\node (r0) [label={0:$\emptyset$}] {} [grow=up]
child {node [label={0:$1$}] {}
}
child {node  [label={180:$0$}] {}
	child {node [label={0:$01$}] {}
	}
	child {node [label={180:$00$}] {}
	}
}
;
\node at (r0) [label={[label distance = 1mm]270:$Q_{2}$}] {};
\end{tikzpicture}
\end{center}
\addvspace{2ex}
\begin{center}
\begin{tikzpicture}
[line width=1pt,level distance=8mm,
every node/.style={draw, color=black, line width=1pt, fill=white, circle, minimum size=2mm, inner sep=0cm, label distance=0mm},
level 1/.style={sibling distance=8mm},
level 2/.style={sibling distance=8mm},
level 3/.style={sibling distance=8mm},
level 4/.style={sibling distance=8mm}]
\node (r0) [label={0:$\emptyset$}] {} [grow=up]
child {node [label={0:$1$}] {}
	child {node [label={0:$11$}] {}
		child {node [label={70:$111$}] {}
		}
		child {node [label={110:$110$}] {}
		}
	}
	child [missing]
}
child {node [label={180:$0$}] {}
}
;
\node at (r0) [label={[label distance = 1mm]270:$Q_{3}$}] {};
\end{tikzpicture}
\hspace{.3in}
\begin{tikzpicture}
[line width=1pt,level distance=8mm,
every node/.style={draw, color=black, line width=1pt,fill=white,circle, minimum size=2mm,  inner sep=0cm, label distance=0mm},
level 1/.style={sibling distance=8mm},
level 2/.style={sibling distance=8mm},
level 3/.style={sibling distance=8mm},
level 4/.style={sibling distance=8mm}]
\node [label={0:$\emptyset$}] {} [grow=up]
child {node [label={0:$1$}] {}
}
child {node [label={180:$0$}] {}
	child [missing]
	child {node [label={180:$00$}] {}
		child {node [label={70:$001$}] {}
		}
		child {node [label={110:$000$}] {}
		}
	}
}
;
\node at (r0) [label={[label distance = 1mm]270:$Q_{4}$}] {};
\end{tikzpicture}
\hspace{.3in}
\begin{tikzpicture}
[line width=1pt,level distance=8mm,
every node/.style={draw, color=black, line width=1pt, fill=white, circle, minimum size=2mm, inner sep=0cm, label distance=0mm},
level 1/.style={sibling distance=8mm},
level 2/.style={sibling distance=8mm},
level 3/.style={sibling distance=8mm},
level 4/.style={sibling distance=8mm}]
\node (r0) [label={0:$\emptyset$}] {} [grow=up]
child {node [label={0:$1$}] {}
	child [missing]
	child {node [label={0:$10$}] {}
		child {node [label={70:$101$}] {}
		}
		child {node [label={110:$100$}] {}
		}
	}
}
child {node [label={180:$0$}] {}
}
;
\node at (r0) [label={[label distance = 1mm]270:$Q_{5}$}] {};
\end{tikzpicture}
\hspace{.3in}
\begin{tikzpicture}
[line width=1pt,level distance=8mm,
every node/.style={draw, color=black, line width=1pt, fill=white, circle, minimum size=2mm, inner sep=0cm, label distance=0mm},
level 1/.style={sibling distance=8mm},
level 2/.style={sibling distance=8mm},
level 3/.style={sibling distance=8mm},
level 4/.style={sibling distance=8mm}]
\node [label={0:$\emptyset$}] {} [grow=up]
child {node [label={0:$1$}] {}
}
child {node [label={180:$0$}] {}
	child {node [label={180:$01$}] {}
		child {node [label={70:$011$}] {}
		}
		child {node [label={110:$010$}] {}
		}
	}
	child [missing]
}
;
\node at (r0) [label={[label distance = 1mm]270:$Q_{6}$}] {};
\end{tikzpicture}
\end{center}
\addvspace{2ex}
\begin{center}
\begin{tikzpicture}
[line width=1pt,level distance=8mm,
every node/.style={draw, color=black, line width=1pt, fill=white, circle, minimum size=2mm, inner sep=0cm, label distance=0mm},
level 1/.style={sibling distance=18mm},
level 2/.style={sibling distance=8mm},
level 3/.style={sibling distance=8mm},
level 4/.style={sibling distance=8mm}]
\node (r0) [label={0:$\emptyset$}] {} [grow=up]
child {node [label={0:$1$}] {}
	child {node [label={90:$11$}] {}
	}
	child {node [label={90:$10$}] {}
	}
}
child {node [label={180:$0$}] {}
	child {node [label={90:$01$}] {}
	}
	child {node [label={90:$00$}] {}
	}
}
;
\node at (r0) [label={[label distance = 1mm]270:$Q_{7}$}] {};
\end{tikzpicture}
\end{center}
\end{center}
\caption{Applying Ramsey with Seven Binary Trees\label{fig:Q's}}
\end{figure}
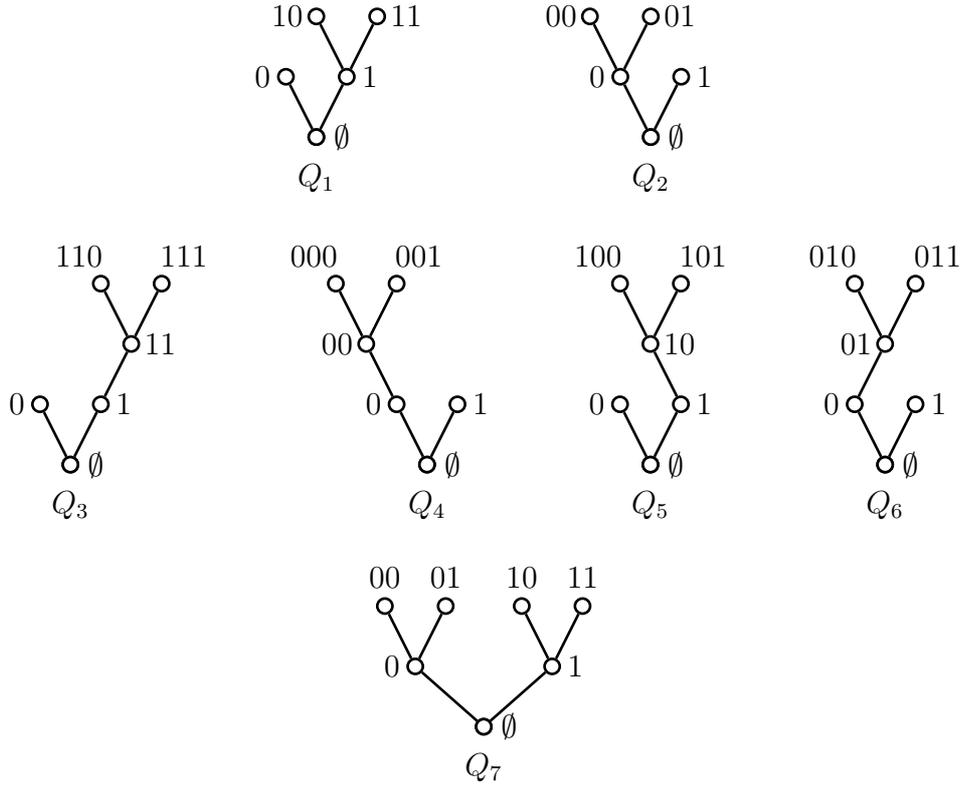

\begin{enumerate}
\item $V_1$ and $V_2$ are adjacent with $z_1=x_2$.
In this case, we associate the pair $(V_1,V_2)$ with a strong
copy of the poset $Q_1$.
\item $V_1$ and $V_2$ are adjacent with $y_1=x_2$.
In this case, we associate the pair $(V_1,V_2)$ with a strong
copy of the poset $Q_2$.
\item $V_1$ and $V_2$ are non-adjacent with $x_2$ in the
right tree above $z_1$.
In this case, we associate the pair $(V_1,V_2)$ with a strong
copy of the poset $Q_3$.
\item $V_1$ and $V_2$ are non-adjacent with $x_2$ in the
left tree above $y_1$.
In this case, we associate the pair $(V_1,V_2)$ with a strong
copy of the poset $Q_4$.
\item $V_1$ and $V_2$ are non-adjacent with $x_2$ in the
left tree above $z_1$.
In this case, we associate the pair $(V_1,V_2)$ with a strong
copy of the poset $Q_5$.
\item $V_1$ and $V_2$ are non-adjacent with $x_2$ in the
right tree above $y_1$.
In this case, we associate the pair $(V_1,V_2)$ with a strong
copy of the poset $Q_6$.
\item $V_1$ and $V_2$ are non-adjacent and there is a vertex
$w$ in $T_n$ so that $x_1$ is in the left tree above
$w$ while $x_2$ is in the right tree above $w$.
In this case, we associate the pair $(V_1,V_2)$ with a strong
copy of the poset $Q_7$.
\end{enumerate}

Also, given a pair $(V_1,V_2)$ of distinct $V$'s in $G_n$,  
there are $6$ ways the intervals $[a_1,b_1]$ and
$[a_2,b_2]$ can appear in the path-decomposition:

\begin{align*}
a_1<a_2<b_1<b_2&\quad\text{Overlapping, moving right}\\
a_2<a_1<b_2<b_1&\quad\text{Overlapping, moving left}\\
a_1<b_1<a_2<b_2&\quad\text{Disjoint, moving right}\\
a_2<b_2<a_1<b_1&\quad\text{Disjoint, moving left}\\
a_1<a_2<b_2<b_1&\quad\text{Inclusion, second in first}\\
a_2<a_1<b_1<b_2&\quad\text{Inclusion, first in second}
\end{align*}

In the arguments to follow, we will abbreviate these $6$ options
as $\OMR$, $\OML$, $\DMR$, $\DML$, $\ISF$ and $\IFS$, respectively.

We then define for each $i\in[7]$ a $6$-coloring $\Phi_i$ of the strong
copies of $Q_i$ in $T_n$.  The colors will be
the six labels $\{\OMR,\OML,\dots,\IFS\}$ listed above.  When
$i\in[7]$ and $Q$ is a strong copy of $Q_i$, then $Q$ is
associated with a pair $(V_1,V_2)$ of vertices from $G_n$.
It is then natural to set $\Phi_i(Q)$ as the label describing
how the pair $([a_1,b_1],[a_2,b_2])$ of intervals are positioned in the
path decomposition.

Now let $p=4\cdot 2^r$.  By iterating on Theorem~\ref{thm:ramsey-bt}, 
we may assume that $n$ is sufficiently large to guarantee 
that there is a subposet $R$ of $T_n$ and
a vector $(\alpha_1,\alpha_2,\dots,\alpha_7)$ of colors such
that $R$ is a strong copy of $T_p$ and for each
$i\in[7]$, $\Phi_i$ assigns color $\alpha_i$ to all strong copies 
of $Q_i$ in $R$.  In the remainder of the argument, we will
abuse notation slightly and simply consider that $R=T_p$.

\noindent
\textbf{Claim 1.}\quad
$\alpha_1$ is either $\OMR$ or $\OML$.

\begin{proof}
A pair $(V_1,V_2)$ of vertices in $G_n$ associated with
a strong copy of $Q_1$ in $T_p$ is adjacent in $G_n$
so that $[a_1,b_1]$ and $[a_2,b_2]$ intersect.  So
$\alpha_1$ cannot be $\DMR$ or $\DML$.
We assume that $\alpha_1$ is $\ISF$ and argue to a contradiction.
The argument when $\alpha_1$ is $\IFS$ is symmetric.
Consider the subposet of $T_p$ consisting of all non-empty strings 
for which each bit, except possibly the last, is a $1$.
We suggest how this subposet appears (at least for a modest value) 
in Figure~\ref{fig:shift}. 

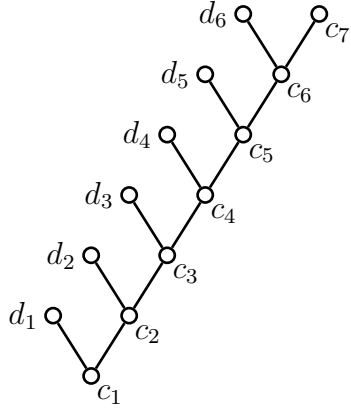
\begin{figure}
\begin{center}
\begin{tikzpicture}
[line width=1pt,level distance=8mm,
every node/.style={draw, color=black, line width=1pt, fill=white, circle, minimum size=2mm, inner sep=0cm, label distance=0mm},
level 1/.style={sibling distance=10mm},
level 2/.style={sibling distance=10mm},
level 3/.style={sibling distance=10mm},
level 4/.style={sibling distance=10mm},
level 5/.style={sibling distance=10mm},
level 6/.style={sibling distance=10mm},
level 7/.style={sibling distance=10mm}
]
\node (r0) [label={-45:$c_1$}] {} [grow=up]
child {node [label={-45:$c_2$}] {}
	child {node [label={-45:$c_3$}] {}
		child {node [label={-45:$c_4$}] {}
			child {node [label={-45:$c_5$}] {}
				child {node [label={-45:$c_6$}] {}
					child {node [label={-45:$c_7$}] {}
					}
					child {node [label={180:$d_6$}] {}
					}
				}
				child {node [label={180:$d_5$}] {}
				}
			}
			child {node [label={180:$d_4$}] {}
			}
		}
		child {node [label={180:$d_3$}] {}
		}
	}
	child {node [label={180:$d_2$}] {}
	}
}
child {node [label={180:$d_1$}] {}
}
;
\end{tikzpicture}
\end{center}
\caption{A Shift Graph in $G_n$\label{fig:shift}}
\end{figure}

Using the labelling given in Figure~\ref{fig:shift},
for each interval $[i,j]$ with
$1\le i<j\le p$, we consider the vertex $V[i,j]=(c_i,d_i,c_j)$.
Clearly, $V[i,j]$ is adjacent to $V[j,k]$ when $1\le i<j<k\le p$,
i.e, these vertices form the shift graph $S_p$.

Let $[a,b]=[a_{V[p-1,p]},b_{V[p-1,p]}]$ be the interval for the 
vertex $V[p-1,p]$.  We claim
that for each $[i,j]$ with $1\le i<j\le p-1$, the interval
for $V[i,j]$ in the path-decomposition for $H_n$ contains $[a,b]$. 
This is immediate if $j=p-1$,
since $(V[i,p-1],V[p-1,p])$ is assigned color $\ISF$.
Now suppose $j<p-1$.  Then $(V[i,j],V[j,p-1])$ is also
$\ISF$, so that in the path-decomposition, the
interval for $V[j,p-1]$ is included in the interval for $V[i,j]$.
By transitivity, we conclude that $[a,b]$ is included 
in the interval for $V[i,j]$.

So the $V$'s in $\{V[i,j]:1\le i<j\le p-1\}$ form a copy
of the shift graph $S_{p-1}$, and all of them are in the
bag $G_n(a)$.  Since $p=4\cdot 2^r$, this is a 
contradiction. 
\end{proof}

Without loss of generality, we take $\alpha_1$ to be $\OMR$,
since if $\alpha_1$ is $\OML$,
we may simply reverse the entire path-decomposition.
To help keep track of the configuration information as
it is discovered, we list this statement as a property.

\noindent
\textbf{Property 1.}\quad $\alpha_1=\OMR$, i.e., $\Phi_1$ assigns
color $\OMR$ to a pair $(V_1,V_2)$ of adjacent vertices in
$G_n$ when $z_1=x_2$.

Although it may not be a surprise, once the color $\alpha_1$ is
set, colors $\alpha_2,\alpha_3,\dots,\alpha_7$ are determined. 

\noindent
\textbf{Property 2.}\quad $\alpha_3=\DMR$,
i.e., $\Phi_3$ assigns color $\DMR$ to
a pair $(V_1,V_2)$ of non-adjacent vertices in $G_n$ when
$x_2$ is in the right tree above $z_1$.

\begin{proof}
Let $(V_1,V_2)$ be a pair of non-adjacent vertices in $G_n$ with
$x_2$ in the right tree above $z_1$.
Then let $w_3$ be the string formed by attaching a $0$ at
the end of $z_1$, and set $V_3=(z_1,w_3,x_2)$. Then $V_3$ is
adjacent to both $V_1$ and $V_2$.  Furthermore,
$\Phi_1(V_1,V_3)=\OMR$ and $\Phi_1(V_3,V_2)=\OMR$.
Accordingly, $\alpha_3$ is either $\OMR$ or $\DMR$. 
We assume that $\alpha_3=\OMR$ and argue to a contradiction.

Consider the shift graph used in the proof of Claim~1.  
Let $a=a_{V[p-1,p]}$ be the left endpoint of the 
interval for $V[p-1,p]$ in the path-decomposition.
We claim that $a$ is in the interval for $V[i,j]$ in the
path-decomposition whenever $1\le i<j\le p-1$.  Again, this 
holds when $j=p-1$ since $\Phi_1(V[i,p-1],V[p-1,p])=\OMR$.  
Also, when $j<p-1$, the color assigned by $\Phi_3$ to the pair 
$(V[i,j],V[p-1,p])$ is also $\OMR$, so that the interval for
$V[i,j]$ in the path-decomposition also contains $a$.  This now implies 
that $G_n(a)$ contains the shift graph $S_{p-1}$.
The contradiction completes the proof.
\end{proof}

\noindent
\textbf{Property 3.}\quad $\alpha_2=\OML$, i.e,
$\Phi_2$ assigns color $\OML$ to a pair $(V_1,V_2)$ of adjacent
vertices in $G_n$ when $y_1=x_2$.  Also, $\alpha_4=\DML$, i.e.,
$\Phi_4$ assigns color $\DML$ to a pair $(V_1,V_2)$ of non-adjacent vertices
in $G_n$ when $x_2$ is in the left tree above $y_1$.

\begin{proof}
We can repeat the arguments given previously to conclude that
one of two cases must hold:\quad Either (1)~$\alpha_2=
\OMR$ and $\alpha_4=\DMR$, or
(2)~$\alpha_2=\OML$ and $\alpha_4=\DML$.
We assume that $\alpha_2=\OMR$ and $\alpha_4=\DMR$ and
argue to a contradiction.  Consider the binary
tree contained in $T_p$ as shown on the left side of 
Figure~\ref{fig:midpoint}.  
Let $V_1=(f,g,h)$, $V_2=(i,j,k)$, 
$V_3=(c,f,e)$ and $V_4=(c,d,i)$.  

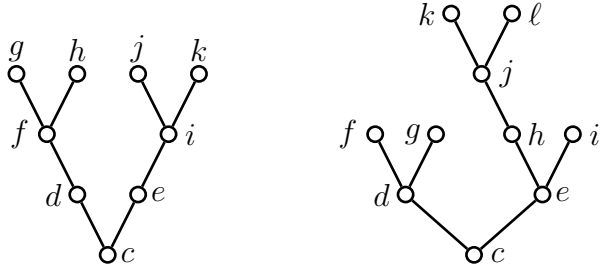
\begin{figure}
\begin{center}
\begin{tikzpicture}
[line width=1pt,level distance=8mm,
every node/.style={draw, color=black, line width=1pt, fill=white, circle, minimum size=2mm, inner sep=0cm, label distance=0mm},
level 1/.style={sibling distance=8mm},
level 2/.style={sibling distance=8mm},
level 3/.style={sibling distance=8mm},
level 4/.style={sibling distance=8mm}]
\node (r0) [label={0:$c$}] {} [grow=up]
child {node [label={0:$e$}] {}
	child {node [label={0:$i$}] {}
		child {node [label={90:$k$}] {}
		}
		child {node [label={90:$j$}] {}
		}
	}
	child [missing]
}
child {node [label={180:$d$}] {}
	child [missing]
	child {node [label={180:$f$}]{}
		child {node [label={90:$h$}] {}
		}
		child {node [label={90:$g$}] {}
		}
	}
}
;
\end{tikzpicture}
\hspace{.5in}
\begin{tikzpicture}
[line width=1pt,level distance=8mm,
every node/.style={draw, color=black, line width=1pt, fill=white, circle, minimum size=2mm, inner sep=0cm, label distance=0mm},
level 1/.style={sibling distance=18mm},
level 2/.style={sibling distance=8mm},
level 3/.style={sibling distance=8mm},
level 4/.style={sibling distance=8mm}]
\node [label={[label distance=.5mm]0:$c$}] {} [grow=up]
child {node [label={0:$e$}] {}
	child {node [label={0:$i$}] {}}
	child {node [label={0:$h$}] {}
		child [missing]
		child {node [label={0:$j$}] {}
			child {node [label={0:$\ell$}] {}
			}
			child {node [label={180:$k$}] {}
			}
		}
	}
}
child {node [label={180:$d$}] {}
	child {node [label={180:$g$}] {}
	}
	child {node [label={180:$f$}] {}
	}
}
;
\end{tikzpicture}
\end{center}
\caption{Two Useful Small Examples\label{fig:midpoint}}
\end{figure}

Since $\Phi_4(V_4,V_1)=\DMR$, we know $b_4<a_1$.  Since
$\Phi_1(V_4,V_2)=\OMR$, we know $a_2<b_4$, so $a_2<a_1$.
Since $\Phi_3(V_3,V_2)=\DMR$, we know $b_3<a_2$ so $b_3<a_1$.
But $\Phi_2(V_3,V_1)=\OMR$, which requires $a_1<b_3$.
The contradiction completes the proof of Property~3.
\end{proof}

\noindent
\textbf{Property 4.}\quad $\alpha_7=\DMR$, i.e., $\Phi_7$ assigns
color $\DMR$ to a pair $(V_1,V_2)$ of non-adjacent vertices in
$G_n$ when there is a vertex $w$ in $T_n$ such that $x_1$ is
in the left tree above $w$ while $x_2$ is in the right tree
above $w$.

\begin{proof}
We again consider the binary tree shown on the left
side of Figure~\ref{fig:midpoint}.  Again, we take
$V_1=(f,g,h)$ and $V_2=(i,j,k)$.  Noting that
$f$ is in the left tree above $c$ and $i$ is in the right
tree above $a$, $\Phi_7(V_1,V_2)=\alpha_7$. 

Now let $V_5=(c,d,e)$. Then $\Phi_4(V_5,V_1)=\DML$ and
$\Phi_3(V_5,V_2)=\DMR$.  These statements imply $\alpha_7=\DMR$.
\end{proof}

\noindent
\textbf{Property 5.}\quad $\alpha_5=\alpha_6=\ISF$, i.e., $\Phi_5$ 
assigns color $\ISF$ to a pair $(V_1,V_2)$ of non-adjacent vertices in
$G_n$ when $x_2$ is in the left tree above $z_1$ and
$\Phi_6$ assigns this pair color $\IFS$ when $x_2$ is in the right
tree above $y_1$.

\begin{proof}
We prove that $\alpha_5=\ISF$.  The argument to show that
$\alpha_6=\ISF$ is symmetric.  Consider the binary tree
shown on the right side of Figure~\ref{fig:midpoint}. 
Let $V_1=(c,d,e)$ and $V_2=(j,k,l)$.
Then $j$ is in the left tree above $e$,
so $\Phi_5(V_1,V_2)=\alpha_5$.  

Now set $V_3=(d,f,g)$ and $V_4=(e,h,i)$.
We observe that $\Phi_2(V_1,V_3)=\OML$, $\Phi_7(V_3,V_2)=\DMR$,
$\Phi_1(V_1,V_4)=\OMR$ and $\Phi_4(V_4,V_2)=\DML$. Together,
these statements imply $\alpha_5=\ISF$.
\end{proof}

Up to this point in the proof, our entire focus has been
on the $V$'s in $G_n$.  We now turn our attention to
properties that the $Y$'s in $H_n$ must satisfy.

\begin{figure}
\begin{center}

\begin{tikzpicture}[every node/.style={draw, color=black, line width=1pt,fill=white,circle, minimum size=2mm,  inner sep=0cm, label distance=0mm},line width=1pt]

\pgfmathsetmacro\L{10}
\pgfmathsetmacro\s{6}

\pgfmathsetmacro\p{2*\s}
\pgfmathsetmacro\pt{\p-2}
\pgfmathsetmacro\ph{\p/2}
\pgfmathsetmacro\sqt{1.7320508075688772}
\pgfmathsetmacro\sqtt{\sqt/2}

\coordinate (T) at (0,\L*\sqtt);
\coordinate (B1) at (-\L/2,0);
\coordinate (B2) at (\L/2,0);


\foreach \i[evaluate=\i as \x using \i/\p] in {0,1,2,...,\p}{
	\coordinate (p\i) at ($(B1)!\x!(T)$);
	\coordinate (q\i) at ($(B2)!\x!(T)$);
}

\foreach \i[count=\l] in {0,2,...,\pt}{
	\pgfmathsetmacro\j{\i+1}
	\pgfmathsetmacro\k{\i+2}
	\coordinate (x\l) at (p\i);
  \coordinate (u\l) at (q\j);
	\draw (x\l) -- (u\l);
 
	\ifnum\i=\pt \else \draw (p\k) -- (q\j);\fi
	
  \coordinate (z\l) at ($(x\l)!0.5!(u\l)$);
}

\foreach \i in {1,2,...,\ph}{
		\pgfmathsetmacro\j{int(\i+1)}

		\node[label={south:$x_{\i}$}] at (x\i) {} [grow=up,level distance=3mm,sibling distance=6mm]
				child [missing]
				child {node[label={west:$y_{\i}$}] {}};
		\node[label={south east:$z_{\i}$}] at (z\i) {};

		\ifnum\i=\ph
				\node at (u\i) {};
		\else
				\node at (u\i) {} [grow'=up,level distance=3mm,sibling distance=6mm]
					child [missing]
					child {node[label={south east:$w_{\j}$}] {}
											[grow'=up,level distance=3mm,sibling distance=4mm]
											child {node {}}
											child {node {}}
									};
		\fi
}

\end{tikzpicture}
\end{center}
\caption{The Final Counter-Example\label{fig:Micek-Y}}
\end{figure}
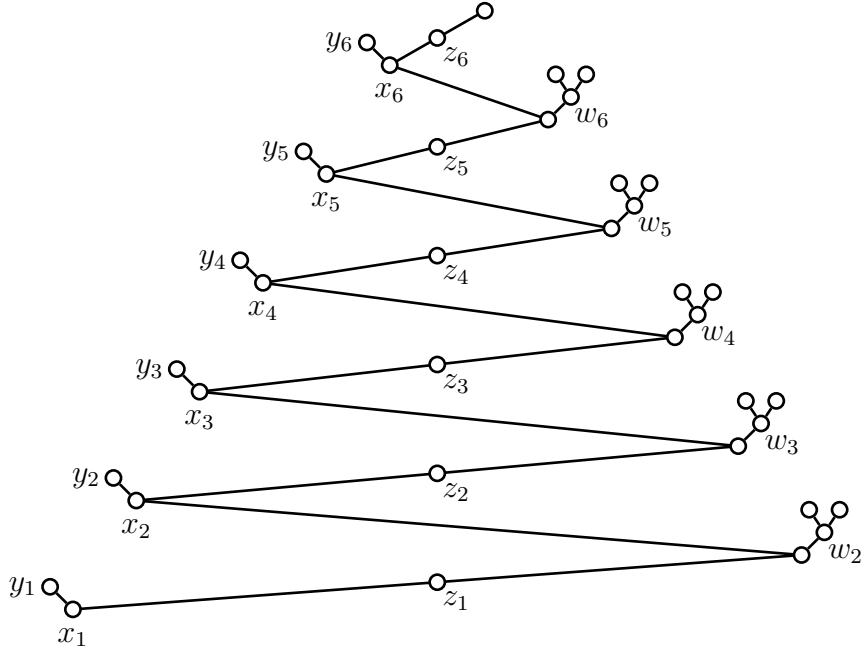

Consider the binary tree shown in Figure~\ref{fig:Micek-Y}.
Of course, we intend that this tree appear inside $T_p$.
In our figure, the ``size'' of this construction is $m=6$, but
since $p=4\cdot 2^r$, we know we can make $m>2^r$.  For each
interval $[i,j]$ with $1\le i<j\le m$, we let $Y[i,j]$ be
the $Y$ whose three leaves are $x_i$, $x_j$ and $w_j$.
Clearly, the family $\{Y[i,j]:1\le i<j\le m\}$ forms a copy
of the shift graph $S_m$.  To reach a final contradiction,
it remains only to show that there
is some integer $k\in\mathbb{N}$ for which
all vertices in $\{Y[i,j]:1\le i<j\le m\}$ belong to $H_n(k)$.

For each $j\in[m]$, we let $V_j=(x_j,y_j,z_j)$, and
as usual, we let $[a_j,b_j]$ be the corresponding interval for
$V_j$ in the path decomposition.  By Property~2, we
have $\alpha_3=\DMR$, so that:

\[
a_1<b_1< a_2<b_2<\dots<a_{m-1}<b_{m-1}<a_m<b_m.
\]

For each $j=2,3,\dots,m$, let $V'_j=(w_j,w_j0,w_j1)$, and
we let $[a'_j,b'_j]$ be the corresponding interval in the 
path-decomposition.  By Property~4, $\alpha_7=\DMR$ so that:

\[
a'_m<b'_m<a'_{m-1}<b'_{m-1}<\dots<a'_3<b'_3<a'_2<b'_2.
\]
Again, since $\alpha_7=\DMR$, we know that $a_m<b_m<a'_m<b'_m$.

Now consider a pair $i,j$ with $1\le i<j\le m$. The vertex
$Y[i,j]$ is adjacent in $H_n$ to both $V_j$ and $V'_j$. This implies
that the interval for $Y[i,j]$ must overlap both $[a_j,b_j]$ and
$[a'_j,b'_j]$.  However, this forces the interval for $Y[i,j]$ to
contain $[b_m,a'_m]$. Therefore, $H_n(b_m)$ contains the
shift graph $S_m$.  With this observation, the proof of
Theorem~\ref{thm:pchrHn} is complete.

We now return to the task of proving Theorem~\ref{thm:pchrGn}, i.e.,
the assertion that $\pchr(G_n)\le 2$ for all $n\ge 1$.  Our
proof for Theorem~\ref{thm:pchrHn} suggests a natural
way to define a path-decomposition of the graph $G_n$ of
$V$'s in the binary tree $T_n$, one that satisfies all five
properties we have developed to this point.  We simply take a
drawing in the plane of $T_n$ using a geometric series
approach.  Taking a standard cartesian coordinate system in
the plane, we place the zero of $T_n$ at the
origin.  If $m\ge 0$ and we have placed a string $x$ of length
$m$ at $(h,v)$, we set $\delta=2^{-m}$ and
place $x1$ and $x0$ at $(h+\delta, v+\delta)$ and
$(h-\delta, v+\delta)$, respectively.

\begin{figure}
\begin{center}
\begin{tikzpicture}[scale=1.4]
\node at (-2.7,2.7) [label={[rectangle,inner sep=0pt]-135:\scriptsize{$V=(x,y,z)$}}] {};
\draw [line width=.5pt, ->, >=stealth] (-5,0)--(5,0); 
\draw [line width=.5pt, ->, >=stealth] (0,0)--(0,4.5); 
\draw [black!30,line width=1mm] (-2.375,3.875)--(-2.25,3.75)--(-3,3)--(-2,2)--(-1,3)--(-1.5,3.5);
\begin{scope}[dashed, line width=.5pt]
\draw (2,0)--(2,2)--(3,2)--(3,3)--(3.5,3)--(3.5,3.5);
\draw (-2.375,0)--(-2.375,3.875);
\draw (-1.5,0)--(-1.5,3.5);
\end{scope}
\draw [line width=2pt] (-1.5,0)--(-2.375,0); 
\begin{scope}[every node/.style={scale=.5}]
\node at (1,0) [above] {1};
\node at (2,1) [right] {1};
\node at (2.5,2) [below] {1/2};
\node at (3,2.5) [right] {1/2};
\node at (3.25,3) [below] {1/4};
\node at (3.5,3.25) [right] {1/4};
\end{scope}
\node at (-2.6,0) [below] {\tiny{$a_V\!=\!\pi(y)$}};
\node at (-1.275,0) [below] {\tiny{$b_V\!=\!\pi(z)$}};
\begin{scope}
[line width=1pt,
every node/.style={draw, color=black, line width=1pt, fill=white, circle, minimum size=1mm, inner sep=0cm, label distance=0mm},
level 1/.style={sibling distance=40mm,level distance=20mm},
level 2/.style={sibling distance=20mm,level distance=10mm},
level 3/.style={sibling distance=10mm,level distance=5mm},
level 4/.style={sibling distance=5mm,level distance=2.5mm},
level 5/.style={sibling distance=2.5mm,level distance=1.25mm}]

\node [label={[rectangle,inner sep=2pt]-90:$\emptyset$}] {} [grow=up]
child {node {}
	child {node {}
		child {node {}
			child {node {}
			}
			child {node {}
			}
		}
		child {node {}
			child {node {}
			}
			child {node {}
			}
		}
	}
	child {node {}
		child {node {}
			child {node {}
			}
			child {node {}
			}
		}
		child {node  {}
			child {node {}
			}
			child {node {}
			}
		}
	}
}
child {node [fill=black, label={[rectangle,fill=white,inner sep=1pt]-135:\tiny{$x\!=\!0$}}] {}
	child {node {}
		child {node {}
			child {node {}
			}
			child {node {}
			}
		}
		child {node [fill=black, label={[rectangle,inner sep=1pt]-135:\tiny{$z\!=\!010$}}] {}
			child {node {}
			}
			child {node {}
			}
		}
	}
	child {node {}
		child {node {}
			child {node {}
				child {node {}}
				child {node[fill=black, label={[rectangle,inner sep=1pt]135:\tiny{$y\!=\!00110\mkern-10mu$}}] {}}
			}
			child {node {}
			}
		}
		child {node {}
			child {node {}
			}
			child {node {}
			}
		}
	}
}
;
\end{scope}
\end{tikzpicture}
\end{center}
\caption{A Path-Decomposition of $G_n$\label{fig:geometry}}
\end{figure}
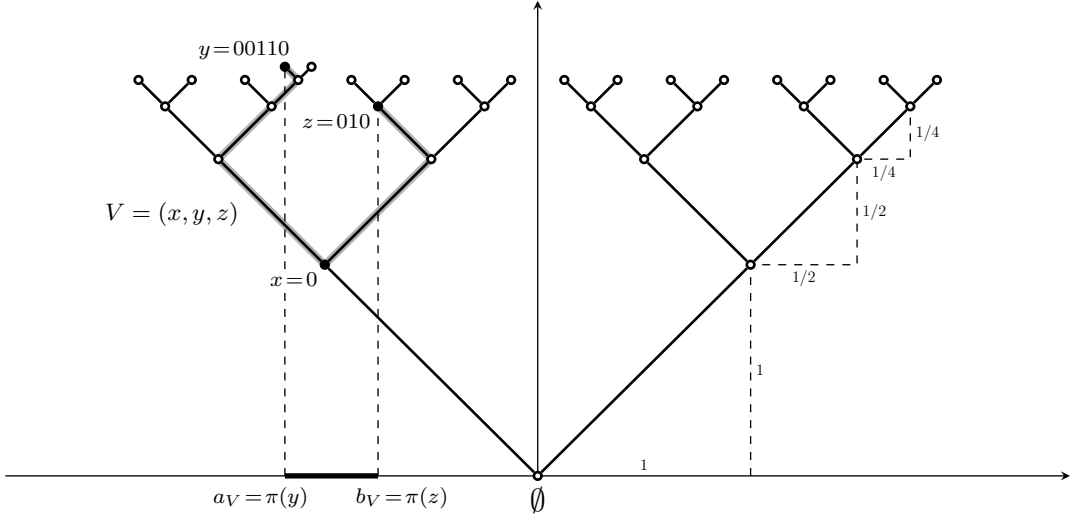

For each $x$ in $T_n$, let
$\pi(x)$ denote the vertical projection of $x$ down onto the horizontal
axis.  In turn, for each $V=(x,y,z)$, we take $a_V=\pi(y)$ and
$b_V=\pi(z)$. To illustrate this construction, we show  
in Figure~\ref{fig:geometry} the interval $[a_V,b_V]$ 
corresponding to the vertex $V=(0,00110,010)$ in $G_n$. 

Clearly, we may consider the host path $P$ for the decomposition as consisting
of all points on the horizontal axis of the form $\pi(x)$ where
$x\in T_n$.  Also, in the natural manner, $\pi(x)$ is adjacent to
$\pi(x')$ in $P$ when there is no string $x''\in T_n$ with
$\pi(x'')$ between $\pi(x)$ and $\pi(x')$.  

So let $x_0\in T_n$ and consider the bag $B = B_{\pi(x_0)}$ consisting 
of all vertices $V=(x,y,z)$ in $G_n$ with $\pi(y)\le \pi(x_0)\le \pi(z)$.
We partition $B$ as $C_1\cup C_2\cup C_3$ where:

\begin{enumerate}
\item A vertex $V=(x,y,z)$ of $B$ belongs to $C_1$ if $\pi(x)<\pi(x_0)$.
\item A vertex $V=(x,y,z)$ of $B$ belongs to $C_2$ if $\pi(x)>\pi(x_0)$.
\item A vertex $V=(x,y,z)$ of $B$ belongs to $C_3$ if $\pi(x)=\pi(x_0)$.
In this case, $x=x_0$.
\end{enumerate}

We now explain why $C_1$, $C_2$
and $C_3$ are independent sets in $G_n$.  This is trivial for
$C_3$.  We give the argument for $C_1$, noting that the argument
for $C_2$ is symmetric.  

Suppose that $V_1$ and $V_2$ are 
adjacent vertices in $C_1$.  If the pair $(V_1,V_2)$ determines
a strong copy of $Q_1$, then $\pi(z_1)=\pi(x_2)<\pi(x_0)$, which
is a contradiction.  On the other hand, if the pair $(V_1,V_2)$ determines
a strong copy of $Q_2$, then $y_1=x_2$ so that
$\pi(y_1)=\pi(x_2)<\pi(x_1)<\pi(x_0)$.  Now the geometric series nature
of the construction implies that $\pi(z_2)<\pi(x_1)<\pi(x_0)$, which is 
again a contradiction.  

With these observations, we have now proved that $\pchr(G_n)\le 3$ for
all $n\ge1$.  This inequality is tight as evidenced by the
following five elements of $G_n$ which form a $5$-cycle:\quad
$V_1=(\emptyset,0,1)$, $V_2=(1,10,11)$, $V_3=(10,100,101)$, $V_4=(101,1010,
1011)$ and $V_5=(1,101,11)$.  Note that $\pi(101)$ is in $[a_i,b_i]$ for
each $i\in[5]$.

Nevertheless, we are able to make a small but important change in the
path-decomposition to obtain a decomposition witnessing that
$\pchr(G_n)\le 2$.  For the integer $n$, let $\epsilon=2^{-2n}$.
Then for each vertex $V=(x,y,z)$ of $G_n$, we change the
interval in the path decomposition for $V$ from $[\pi(y),\pi(z)]$ to
$[\pi(y)+\epsilon,\pi(z)-\epsilon]$.  Our choice of $\epsilon$ guarantees
that we still have a path-decomposition of $G_n$.

Again, we consider an element $x_0$ of $T_n$ and the bag $B$ consisting of
all $V=(x,y,z)$ with $\pi(y)+\epsilon\le \pi(x_0)\le \pi(z)-\epsilon$.  
As before, $C_1$, $C_2$ and $C_3$ are independent sets, although membership in
these three sets has been affected by the revised path-decomposition.
We claim that $C_1\cup C_3$ is also an independent set, so that the partition
$B=(C_1\cup C_3)\cup C_2$ witnesses that $\pchr(G_n)\le 2$.

Suppose to the contrary that 
$V_1\in C_1$ and $V_3\in C_3$ with
$V_1$ adjacent to $V_3$ in $G_n$. Clearly, this requires that
$(V_1,V_3)$ is associated with a strong copy of the binary tree
$Q_1$ as shown in Figure~\ref{fig:Q's}.  This implies that $z_1=x_0=x_3$
so that $b_1=\pi(z_1)-\epsilon=\pi(x_0)-\epsilon$, which contradicts the
assumption that $a_1<\pi(x_1)<\pi(x_0)\le b_1$. 
The contradiction completes the proof of Theorem~\ref{thm:pchrGn}.

\section{Acknowledgements}
We thank an anonymous editor for pointing out the Milliken's result in~\cite{bib:Mill}.
In an earlier version of the paper we were proving Theorem~\ref{thm:ramsey-bt} on our own.

\end{document}